\documentclass{elsart}

\usepackage{natbib}

\usepackage{amssymb}

\begin{document}

\begin{frontmatter}



\title{Cross Validation for Comparing Multiple Density Estimation Procedures}


\author{Heng Lian}
\ead{Heng\_Lian@brown.edu}
\address{Division of Applied Mathematics\\Brown University\\182 George St\\Providence, RI 02912}

\begin{abstract}
We demonstrate the consistency of cross validation for comparing multiple density estimators using simple inequalities on the likelihood ratio. In nonparametric problems, the splitting of data does not require the domination of test data over the training/estimation data, contrary to Shao (1993). The result is complementary to that of Yang (2005) and Yang (2006).
\end{abstract}

\begin{keyword}
Density estimation\sep Cross validation\sep Kullback-Leibler divergence\sep Hellinger distance

\end{keyword}

\end{frontmatter}

\section{Introduction}
\label{}
Cross validation (CV) is a common procedure used for smoothing parameter estimation or model selection. When multiple procedures are being compared, it is customary to choose the one that obtains the smallest ``loss", which is defined specifically to the problem at hand. If both the procedure used to obtain the estimate and the computation of the loss are based on the same set of data, it is a well-known effect that the estimated loss is biased due to the double use of the same observations. In order to obtain an unbiased estimate, one approach is to use a penalty term that takes into account the complexity of the model. This approach includes AIC, BIC, $C_p$, etc. 
A simpler approach, which is closely related, when we have the luxury of enough observations,   is to split the data in such way that one part is used to obtain the estimate, and the other separate hold-out data is used to evaluate the loss. The main advantage of this approach is that it can be easily applied in various situations (without theoretical derivation) to select one out of many competing procedures. There are a few different variations to this approach, including leave-one-out CV, k-fold CV, Generalized CV, etc. A  well-studied case is the linear regression problem, where Shao (1993) showed the surprising result that leave-one-out CV will select the models with extra redundant variables with nonvanishing probability
(This theory assumed that the number of covariates is fixed. It is a different story when the number of covariates grows with $n$). In order to restore consistency, one should split the data such that the size of the evaluation part of the data is dominating. All of the above techniques are summarized and compared in Shao (1997) in the context of linear regression with different kinds of asymptotics.

Yang (2005) studied the problem of cross-validation in the context of nonparametric regression comparing a finite number of estimators. It is shown in that paper that under the $L_2$ loss, as long as one of the competing procedures converges at a nonparametric rate, the dominance of evaluation data is not necessary for consistency. Instead, the two parts of the data can be of the same order, which is surprising considering the corresponding result for linear regression.
The proof of Yang's result is based on an application of Bernstein's inequality. Similarly, it is shown in Yang (2006) that cross validation is consistent in classification problems where the consistency also depends on the rate of disagreement between the two classifiers.

In this paper, we consider the problem of density estimation when the observations are generated i.i.d. from the true distribution $P_0$. There exists a large literature on density estimators, earlier results focus on linear estimators including kernel density estimator, later developments include wavelets thresholding and adaptive width kernel that achieve minimax rate of convergence in a large class of Besov spaces where no linear estimate can attain the optimal convergence rate. Faced with such large choices of estimator with different theoretical properties, it is important to select one that has optimal performance for the current problem. Cross validation can be directly applied by splitting the observations into two groups. Different estimators, such as kernel estimates and wavelets, can be obtained based on the first part of data, then the likelihood of the second group of data can be evaluated and compared. Finally, the estimator that obtains the largest likelihood is chosen as the winner. A natural question is whether this process will return the optimal procedure. In particular, what condition on the splitting ratio should be satisfied in order to ensure the consistency property? 
  The main conclusion in this paper is similar to that of Yang (2005), that is, we do not necessarily need to assume the size dominance of the evaluation data as in linear regression problem. 

\section{Consistency result}
Consider the situation where we have $n$ observations $X_1, X_2,\ldots, X_n$ generated i.i.d. from an underlying distribution $P_0$ with corresponding density $p_0$. There exists many well-known density estimators. For parametric procedures, such as mixture modeling using parametric families of densities, the convergence rate is usually $1/\sqrt{n}$. On the other hand, for nonparametric procedures, the rate of convergence is slower, depending on the smoothness property of the true density. Density estimation is closely related to regression problem, as demonstrated in a series of well-known papers (see, e.g., Brown and Low (1996), Nussbaum (1996)). 

With many possible choices for density estimation procedures, both parametric and nonparametric, one needs to find the best estimator among them for the current data. Parametric procedures have a faster rate of convergence when the model is correct, but suffer from a nonvanishing bias when the true distribution lies outside of the parametric family. Nonparametric procedures are more flexible but lose in efficiency when the underlying density is of a known parametric form. In practice, we need to know which procedure is best without knowledge of the true distribution.

We start by splitting the data into two parts: the estimation data $X^1=\{X_1,\ldots, X_{n_1}\}$ and the evaluation data $X^2=\{X_{n_1+1},\ldots, X_n\}$, and let $n_2=n-n_1$. We assume we have many estimation procedures $\{\hat{P_i}\}_{i=1}^{m_n}$ (note the number of potential choices can grow with $n$), which will produce density estimates $\{\hat{p}_i^{(n_1)}(x;X_1,\ldots,X_{n_1})\}$, we will omit the dependence of $\hat{p}_i^{(n_1)}$ on the training data $X^1$ in the following. To choose the best procedure among those $m_n$, the test data $X^2$ is used to evaluate the likelihood: $\hat{p}_i^{(n_1)}(X^2)=\prod_{k=n_1+1}^{n}\hat{p}_i^{(n_1)}(X_k)$. If $\hat{p}_j^{(n_1)}(X^2)=\max_i\hat{p}_i^{(n_1)}(X^2)$, then the procedure $\hat{P}_j^{(n_1)}$ is selected as the final estimator. The desired property is that this cross-validation procedure will select the best one with high probability. We will use a loss function $d(p_0,p)$ to measure the closeness of $p$ to the true density $p_0$. In this paper, we will adopt the commonly used Hellinger distance as the loss function: $d_H(p_0,p)=(\int (\sqrt{p_0}-\sqrt{p})^2)^{1/2}$. Another commonly used measure of loss in the context of density estimation is the Kullback-Leibler divergence (which is not a true distance) $d_K(p_0,p)=\int p_0\log\frac{p_0}{p}$. It is always true that $d_H^2\le d_K$. Under some mild assumptions, these two measures of loss are almost equivalent. The simplest case under which this is true is when the class of densities considered are uniformly bounded away from zero and infinity, so that the ratio $\frac{p_0}{p}$ is uniformly bounded, then $d_K(p_0,p)=O(d_H^2(p_0,p))$ (see, e.g., Lemma 8.2 in Ghosal et al. (2000)). The more complicated techniques similar to those used in section 3 of that paper can also be used when the estimate is constrained to be within a finite approximation set (this will result in an extra logarithmic factor). The paper of Yang and Barron (1999) contains more information regarding the relationships between $d_H^2$, $d_K$ and the $L_p$ loss function, and established some equivalence result between them under some conditions. In the rest of the paper, we will assume that $d_K\le M d_H^2$ for some constant $M$. 

\begin{thm}
$P_0(\hat{p}_1^{(n_1)}(X^2)>\hat{p}_i^{(n_1)}(X^2),\forall i>1|X^1)\rightarrow 1$ if the following conditions hold:

(1) $n_1\rightarrow\infty, n_2\rightarrow\infty$

(2) $n_2\min_{i>1}v_{n_1,i}^2\rightarrow\infty,$ and $\log m_n=o(n_2\min_{i>1}v_{n_1,i}^2)$

(3) There exists $c<1$ and $t_n>0$ s.t. $nt_n\rightarrow\infty$ and $\frac{Mv_{n_1,1}^2+s_{n_1}t_{n_2}}{cv_{n_1,i}^2}<1,   \forall i>1$

where $v_{n,i}=d_H(p_0,\hat{p}_i^{(n)})$, $s_n=V(p_0,\hat{p}_1^{(n)})$, and $V(p_0,p)=\int p_0(\log\frac{p_0}{p})^2$
\end{thm}
\begin{rem}
In the statement of the theorem, the probability is conditioned on $X^1$, and $v_{n_1,i}$ and $s_{n_1}$ in  conditions (2) and (3) are random variables, so the theorem should be interpreted as $P_0(\hat{p}_1^{(n_1)}(X^2)>\hat{p}_i^{(n_1)}(X^2),\forall i>1|X^1)\rightarrow 1$ on the set that the conditions (1)-(3) hold.
\end{rem}
Proof. The proof is based on simple likelihood ratio inequalities in Wong and Shen (1995). All the probabilities below are implicitly conditioned on the training data $X^1$. From Lemma 1 in Wong and Shen (1995), we have, for $i\neq 1, b>0$,
\begin{eqnarray*}
 P_0(\frac{\hat{p}_i^{(n_1)}(X^2)}{p_0(X^2)}\ge exp(-n_2b))&\le& exp(\frac{n_2b}{2}-\frac{n_2v_{n_1,i}^2}{2})\\
\end{eqnarray*}
Choosing $b=cv_{n_1,i}^2$ ($c$ as in the above assumption (3)), we get 
\begin{eqnarray*}
P_0(\frac{\hat{p}_i^{(n_1)}(X^2)}{p_0(X^2)}\ge exp(-n_2cv_{n_1,i}^2))&\le& exp(-\frac{n_2v_{n_1,i}^2}{2}(1-c))
\end{eqnarray*}
Denote by $W_{n_2}$ the event $\{\frac{d_K^{(n_2)}(p_0,\hat{p}_1^{(n_1)})-d_K(p_0,\hat{p}_1^{(n_1)})}{V(p_0,\hat{p}_1^{(n_1)})}\ge t_{n_2}\}$, where $d_K^{(n_2)}(p_0,\hat{p}_1^{(n_1)})$ is the empirical version of $d_K(p_0,\hat{p}_1^{(n_1)})$ on evaluation data $X^2$:
\begin{eqnarray*}
d_K^{(n_2)}(p_0,\hat{p}_1^{(n_1)})=\frac{1}{n_2}\sum_{i=n_1+1}^n \log \frac{p_0(X_i)}{\hat{p}_1^{(n_1)}(X_i)}
\end{eqnarray*}
. By Chebyshev's inequality, $P_0(W_{n_2})\le\frac{1}{n_2t_{n_2}}$. 
Denoting $d=Mv_{n_1,1}^2+s_{n_1}t_{n_2}$,
\begin{eqnarray*}
&& P_0(\frac{\hat{p}_1^{(n_1)}(X^2)}{p_0(X^2)}\le exp(-n_2d))\\
&=& P_0(exp\{-n_2 V(p_0,\hat{p}_1^{(n_1)}) \frac{d_K^{(n_2)}(p_0,\hat{p}_1^{(n_1)})-d_K(p_0,\hat{p}_1^{(n_1)})}{V(p_0,\hat{p}_1^{(n_1)})}\}  exp\{-n_2d_K(p_0,\hat{p}_1^{(n_1)})\}\\
&&\;\;\le exp\{-n_2d\})\\
&\le&P_0\left(exp\{-n_2 V(p_0,\hat{p}_1^{(n_1)}) \frac{d_K^{(n_2)}(p_0,\hat{p}_1^{(n_1)})-d_K(p_0,\hat{p}_1^{(n_1)})}{V(p_0,\hat{p}_1^{(n_1)})} \}\le exp\{-n_2(d-Mv_{n_1,1}^2)\}\right)\\
&\le&P_0(W_{n_2})\rightarrow 0
\end{eqnarray*}
The last inequality above holds since on the set
\begin{eqnarray*}
W_{n_2}^c=\{\frac{d_K^{(n_2)}(p_0,\hat{p}_1^{(n_1)})-d_K(p_0,\hat{p}_1^{(n_1)})}{V(p_0,\hat{p}_1^{(n_1)})}<t_{n_2}\},
\end{eqnarray*}
we have 
\begin{eqnarray*}
exp\{-n_2 V(p_0,\hat{p}_1^{(n_1)}) \frac{d_K^{(n_2)}(p_0,\hat{p}_1^{(n_1)})-d_K(p_0,\hat{p}_1^{(n_1)})}{V(p_0,\hat{p}_1^{(n_1)})} \}&>&exp\{-n_2s_{n_1}t_{n_2}\}\\
&=&exp\{-n_2(d-Mv_{n_1,1}^2)\}
\end{eqnarray*}
Note that when $d=Mv_{n_1,1}^2+s_{n_1}t_{n_2}<cv_{n_2,i}^2$ (assumption (3)), $ \frac{\hat{p}_1^{(n_1)}(X^2)}{p_0(X^2)}>e^{-n_2d}$ and $\frac{\hat{p}_i^{(n_1)}(X^2)}{p_0(X^2)}<e^{-n_2cv_{n_2,i}^2}$ implies $\hat{p}_1^{(n_1)}(X^2)>\hat{p}_i^{(n_1)}(X^2)$. So we can bound the probability that the cross validation procedure chooses an estimator other than $\hat{P}_1$:
\begin{eqnarray*}
&&P_0(\hat{p}_1^{(n_1)}(X^2)<\hat{p}_i^{(n_1)}(X^2) \mbox{ for some } i)\\
&\le& P_0(\frac{\hat{p}_1^{(n_1)}(X^2)}{p_0(X^2)}<e^{-n_2d})+\sum_{i=2}^{m_n}P_0(\frac{\hat{p}_i^{(n_1)}(X^2)}{p_0(X^2)}>e^{-n_2cv_{n_2,i}^2})\\
&\le&P_0(W_{n_2})+m_nexp(-\frac{n_2(1-c)}{2}\min_i v_{n_1,i}^2)
\end{eqnarray*}
The above expression converges to zero under the assumed condition (2). $\Box$

\begin{rem}
The above theorem cannot be directly applied since $v_{n_1,i}$ are random variables depending on $X^1$. We will specialize the result to the two procedures case below in Corollary 1.
\end{rem}
\begin{rem}
Under some mild conditions, we will have $s_{n_1}=V(p_0,\hat{p}_1^{(n_1)})=O(d_H^2(p_0,\hat{p}_1^{(n_1)}))$, see, e.g., Theorem 5 in Wong and Shen (1995).
\end{rem}
\begin{rem}
As in Yang (2005), we can consider the case where we have multiple different splittings of the original data. Cross validation as stated above can be applied to each splitting separately, and then use majority vote to choose the final procedure.
\end{rem}
Now we give a definition comparing two procedures by their rate of convergence. 
\begin{defn}
Procedure $\hat{P}_1$ is asymptotically better than $\hat{P}_2$ under the loss function $d_H$ if for some sequence $\epsilon_n\rightarrow 0$
\begin{displaymath}
  \lim_{n\rightarrow \infty}P_0(d_H(p_0,\hat{p}_1^{(n)})<\epsilon_n d_H(p_0,\hat{p}_2^{(n)}))\rightarrow 1
\end{displaymath}
\end{defn}

Under this definition we can state the following corollary, the definition of exact rate of convergence is similar to Definition 3 in Yang (2005).
\begin{cor}
Considering two procedures for density estimation where one is asymptotically better than the other. Suppose the exact rate of convergence of $d_H(p_0,\hat{p}_1^{(n)})$ and $d_H(p_0,\hat{p}_2^{(n)})$  are $p_n$ and $q_n$ respectively. Assume that $V(p_0,\hat{p}_i^{(n_1)})=O(d_K(p_0,\hat{p}_i^{(n_1)})), i=1,2$. If $n_1\rightarrow\infty,n_2\rightarrow\infty, \sqrt{n_2}\max(p_{n_1},q_{n_1})\rightarrow\infty$, then the cross validation is consistent in the sense of choosing the asymptotically better procedure with probability tending to 1.
\end{cor}

\begin{rem}
If one of the procedures has nonparametric rate $n^{-\alpha}$ with $\alpha<1/2$, then $\frac{n_1}{n_2}=O(1)$ will suffice for $\sqrt{n_2}\max(p_{n_1},q_{n_1})\rightarrow\infty$
\end{rem}

\section{Conclusions}
We give a simple proof of the consistency of cross validation in the context of density estimation. Although it is shown in Shao (1993) that leave-one-out cross validation is inconsistent for linear regression problem, it is unclear to us whether this is the case for nonparametric problems. Another interesting problem is that when multiple splittings are  available, we can either use majority voting as in Yang (2005) or choose the procedure with the largest product of individual likelihood for each splitting. The comparison of these two approaches is similar to the tradeoff between model selection and model averaging.


\makeatletter 
\def\@biblabel#1{\hspace*{-\labelsep}} 
\makeatother

\end{document}